\documentclass{article}
\usepackage{amssymb,amsfonts,amsmath,amsthm,enumerate}
\usepackage[numeric]{amsrefs}
\usepackage{mathtools}
\usepackage[nocompress]{cite}

\def \C{{\mathbb C}}

\def \f{\frac}
\def \t{\theta}
\def \z{\zeta}

\def \vp{\varphi}

\def \Q{{\mathbb Q}}

\def \R{{\mathbb R}}
\def \Re{\operatorname{Re}}

\def \Z{{\mathbb Z}}
\def \l{\left}
\def \r{\right}

\def \z{\zeta}

\newtheorem{theorem}{Theorem}
\newtheorem{thm}{Theorem}

\numberwithin{equation}{section} \theoremstyle{definition}

\newtheorem{example}{Example}

\begin{document}

\title{Variations of Ramanujan's Euler Products}
\author{Shin-ya Koyama\footnote{Department of Biomedical Engineering, Toyo University,
2100 Kujirai, Kawagoe, Saitama, 350-8585, Japan.} \ \& Nobushige Kurokawa\footnote{Department of Mathematics, Tokyo Institute of Technology, 
Oh-okayama, Meguro-ku, Tokyo, 152-8551, Japan.}}

\maketitle

\begin{abstract}
We study the meromorphy of various Euler products of degree two attached to 
holomorphic Hecke eigen cusp forms for the elliptic modular group,
including Ramanujan's $\Delta$-function.
\end{abstract}



\section*{Introduction}
In 1916 Ramanujan \cite{R} studied the Euler product associated to Ramanujan's $\tau$-function defined as
$$
\Delta(z)=e^{2\pi i z}\prod_{n=1}^\infty(1-e^{2\pi inz})^{24}
=\sum_{n=1}^\infty\tau(n) e^{2\pi inz}
$$
for $\Im(z)>0$. This $\Delta(z)$ is a Hecke-eigen cusp form of weight 12 for the elliptic modular group $SL(2,\Z)$.
Ramanujan's Euler product is written in the normalized form as
$$
L_\Delta(s)=\prod_{p:\,\text{prime}}(1-a(p)p^{-s}+p^{-2s})^{-1},
$$
where $a(n)=\tau(n)n^{-\f{11}2}.$

Ramanujan conjectured that $|a(p)|\le2$ for all primes $p$ and it was proved by Deligne \cite{D} in 1974.
Hence there exists a unique $\t(p)\in[0,\pi]$ satisfying $a(p)=2\cos(\t(p))$.

Thus
\begin{align*}
L_\Delta(s)
&=\prod_{p:\,\text{prime}}(1-2\cos(\t(p))p^{-s}+p^{-2s})^{-1}\\
&=\prod_{p:\,\text{prime}}\det(1-\alpha_\Delta^0(p)p^{-s})^{-1},
\end{align*}
where
$$
\alpha_\Delta^0(p)
=\l[\begin{pmatrix}
e^{i\t(p)}&0\\
0&e^{-i\t(p)}
\end{pmatrix}\r]
\in\mathrm{Conj}(SU(2)).
$$

More general Euler products were constructed by Serre \cite{Se}, Langlands \cite{L} and Tate \cite{T} as
$$
L_\Delta(s,\mathrm{Sym}^m)
=\prod_{p:\,\text{prime}}\det(1-\mathrm{Sym}^m(\alpha_\Delta^0(p))p^{-s})^{-1}
$$
for the irreducible representation
$\mathrm{Sym}^m:\ SU(2)\to U(m+1)$ with $m=0,1,2,3,...$.
Quite recently, Newton-Thorne \cite{NT1, NT2} proved that there exists an automorphic representation
$\pi_m$ of $GL(m+1,\mathbb A_\Q)$ satisfying
$$
L(s,\pi_m)=L_\Delta(s,\mathrm{Sym}^m)
$$
for each $m=0,1,2,3,...$ and that $\pi_m$ $(m\ge 1)$ are cuspidal corresponding to the holomorphy of
$L_\Delta(s,\mathrm{Sym}^m)$.

The situation is quite similar for each holomorphic Hecke eigen cusp form $\vp$ of general weight for the elliptic modular group.

In this paper we first study the Euler product
$$
Z_m^{\pm}(s)=\prod_{p:\,\text{prime}}(1\pm2\cos(m\t(p))p^{-s}+p^{-2s})^{-1}
$$
for $m=0,1,2,3,...$ and show the following result.

\begin{theorem}
$Z_m^{\pm}(s)$ has an analytic continuation to all $s\in\C$ as a meromorphic function.
\end{theorem}

In fact, this is easily seen from Newton-Thorne \cite{NT1, NT2} since each $Z_m^\pm(s)$ can be written
explicitly in terms of $L_\vp(s,\mathrm{Sym}^n)$.
In particular $Z_m^-(s)$ has a simple expression as
\begin{align*}
Z_m^-(s)=\begin{cases}
\dfrac{L_\vp(s,\mathrm{Sym}^m)}{L_\vp(s,\mathrm{Sym}^{m-2})}&(m\ge2)\\
L_\vp(s)&(m=1).
\end{cases}
\end{align*}

Next, we show the converse in some sense.
Actually we determine all (non-constant) monic polynomials $f(x)\in\Z[x]$ such that
\begin{align*}
Z^\pm(s,f)
&=\prod_{p:\,\text{prime}}(1\pm f(a(p))p^{-s}+p^{-2s})^{-1}\\
&=\prod_{p:\,\text{prime}}(1\pm f(2\cos(\t(p))p^{-s}+p^{-2s})^{-1}
\end{align*}
are meromorphic on $\C$.

\begin{theorem}
Let $f(x)\in\Z[x]$ be a monic polynomial of degree $m\ge1$.
Then the following properties are equivalent.
\begin{enumerate}[\rm(1)]
\item \vskip -5mm
$Z^\pm(s,f)$ is meromorphic on $\C$.
\item $|f(x)|\le2$ for $-2\le x\le 2$.
\item $f(x)=2T_m(\f x2)$. Here, $T_m(x)$ is the Chebyshev polynomial defined by $T_m(\cos\t)=\cos(m\t)$.
\item $Z^\pm(s,f)=Z_m^\pm(s)$.
\end{enumerate}
\end{theorem}

\begin{example}
Let
$$
Z^m(s)=\prod_{p:\,\text{prime}}(1-(a(p^2)-m)p^{-s}+p^{-2s})^{-1}
$$
for $m\in\Z$. Then,
\begin{center}
$Z^m(s)$ is meromorphic on $\C$ \quad$\Longleftrightarrow$\quad $m=1$.
\end{center}
Moreover,
\begin{align*}
L_\vp(s,\mathrm{Sym}^2)
&=\z(s) \prod_{p:\,\text{prime}}(1-(a(p^2)-1)p^{-s}+p^{-2s})^{-1}\\
&=\prod_{p:\,\text{prime}}(1-a(p^2)p^{-s}+a(p^2)p^{-2s}-p^{-3s})^{-1}.
\end{align*}
This expression was used by Shimura \cite{S} for the proof of the holomorphy of $L_\vp(s,\mathrm{Sym}^2)$.
\end{example}

\begin{example}
For $m\in\Z$,
\begin{center}
$\prod\limits_{p:\,\text{prime}}(1-(a(p)-m)p^{-s}+p^{-2s})^{-1}$ is meromorphic on $\C$ \quad$\Longleftrightarrow$\quad $m=0$.
\end{center}
\end{example}

Otherwise, each $Z^\pm(s,f)$ has the natural boundary $\Re(s)=0$ as in the following result.

\begin{theorem}
Let $f(x)\in\Z[x]$ be a non-constant monic polynomial.
Then the following properties are equivalent.
\begin{enumerate}[\rm(1)]
\item \vskip -4mm
$Z^\pm(s,f)$ is not meromorphic on $\C$.
\item $Z^\pm(s,f)$ has an analytic continuation in $\Re(s)>0$ with the natural boundary $\Re(s)=0$.
\end{enumerate}
\end{theorem}

Our results come from a general meromorphy theorem for Euler products shown in \cite{K1, K2} assuming 
\cite{NT1, NT2}. We notice that the proof is similar to the case of meromorphy of Dirichlet series
$\sum\limits_{n=1}^\infty a(n)^m n^{-s}$ and
$\sum\limits_{n=1}^\infty a(n^m) n^{-s}$ treated in \cite{K1, K2}.
For $m=1,2$ they are meromorphic on $\C$ and for $m\ge3$ they are meromorphic in $\Re(s)>0$ with the
natural boundary $\Re(s)=0$. We refer to \cite{KK} for an application to rigidity of Euler products.

\section{Proof of Theorem A}
\subsection{The case $m=0$}
From
$$
Z_0^\pm(s)=\prod_{p:\,\text{prime}}(1\pm2p^{-s}+p^{-2s})^{-1}
=\prod_{p:\,\text{prime}}(1\pm p^{-s})^{-2}
$$
we see that
$$
Z_0^-(s)=\prod_{p:\,\text{prime}}(1- p^{-s})^{-2}=\z(s)^{2}
$$
and
$$
Z_0^+(s)=\prod_{p:\,\text{prime}}(1+ p^{-s})^{-2}
=\prod_{p:\,\text{prime}}\l(\f{1- p^{-2s}}{1- p^{-s}}\r)^{-2}
=\f{\z(2s)^2}{\z(s)^{2}}
$$
are meromorphic on $\C$.

\subsection{The case $m=1$}
By identifying 
$$
Z_1^-(s)=\prod_{p:\,\text{prime}}(1-2\cos(\t(p)) p^{-s}+p^{-2s})^{-1}=L_\vp(s)
$$
and
\begin{align*}
Z_1^+(s)
&=\prod_{p:\,\text{prime}}(1+2\cos(\t(p)) p^{-s}+p^{-2s})^{-1}\\
&=\prod_{p:\,\text{prime}}[(1+e^{i\t(p)}p^{-s})(1+e^{-i\t(p)}p^{-s})]^{-1}\\
&=\prod_{p:\,\text{prime}}\l[\f{(1-e^{2i\t(p)}p^{-2s})(1-e^{-2i\t(p)}p^{-2s})}{(1-e^{i\t(p)}p^{-s})(1-e^{-i\t(p)}p^{-s})}\r]^{-1}\\
&=\f{L_\vp(2s,\mathrm{Sym}^2)}{L_\vp(s,\mathrm{Sym}^1)L_\vp(2s,\mathrm{Sym}^0)},
\end{align*}
where
\begin{align*}
L_\vp(s,\mathrm{Sym}^1)&=L_\vp(s),\\
L_\vp(s,\mathrm{Sym}^0)&=\z(s),
\end{align*}
we know that both are meromorphic on $\C$.

\subsection{The case $m\ge2$}
First
\begin{align*}
Z_m^-(s)
&=\prod_{p:\,\text{prime}}[(1-e^{im\t(p)}p^{-s})(1-e^{-im\t(p)}p^{-s})]^{-1}\\
&=\prod_{p:\,\text{prime}}\l[\f{(1-e^{im\t(p)}p^{-s})(1-e^{i(m-2)\t(p)}p^{-s})\cdots(1-e^{-im\t(p)}p^{-s})}
{(1-e^{i(m-2)\t(p)}p^{-s})\cdots (1-e^{-i(m-2)\t(p)}p^{-s})}\r]^{-1}\\
&=\f{L_\vp(s,\mathrm{Sym}^m)}{L_\vp(s,\mathrm{Sym}^{m-2})}
\end{align*}
is a meromorphic function on $\C$ by means of \cite{BLGHT, NT1, NT2}.

Secondly
\begin{align*}
Z_m^+(s)
&=\prod_{p:\,\text{prime}}[(1+e^{im\t(p)}p^{-s})(1+e^{-im\t(p)}p^{-s})]^{-1}\\
&=\prod_{p:\,\text{prime}}\l[\f
{(1-e^{2im\t(p)}p^{-2s})(1-e^{-2im\t(p)}p^{-2s})}
{(1-e^{im\t(p)}p^{-s})(1-e^{-im\t(p)}p^{-s})}\r]^{-1}\\
&=
\f{L_\vp(2s,\mathrm{Sym}^{2m})}{L_\vp(2s,\mathrm{Sym}^{2m-2})}
\f{L_\vp(s,\mathrm{Sym}^{m-2})}{L_\vp(s,\mathrm{Sym}^{m})}
\end{align*}
is meromorphic on $\C$ using \cite{BLGHT, NT1, NT2}.
\hfill\qed

\section{Meromorphy of Euler products}
We recall the meromorphy result of \cite[\S3 Theorem 7]{K2} needed in this paper.
Take the triple
$$
E=(P(\Z),\ SU(2),\ \alpha_\vp^0),
$$
where $P(\Z)$ is the set of prime numbers, and the map
$$
\alpha_\vp^0:\ P(\Z)\to\mathrm{Conj}(SU(2))
$$
defined as
$$
\alpha_\vp^0(p)=\l[\begin{pmatrix}
e^{i\t(p)}&0\\
0&e^{-i\t(p)}\end{pmatrix}\r].
$$
We identify $[0,\pi]$ and $\mathrm{Conj}(SU(2))$ via the following map
$$
\begin{matrix}
[0,\pi] & \longrightarrow & \mathrm{Conj}(SU(2)).\\
\rotatebox{90}{$\in$} & &\rotatebox{90}{$\in$}\\
\t & \longmapsto & \l[\begin{pmatrix}
e^{i\t}&0\\
0&e^{-i\t}\end{pmatrix}\r]
\end{matrix}
$$
Notice that the normalized measure on $\mathrm{Conj}(SU(2))=[0,\pi]$ coming from the normalized Haar measure of $SU(2)$ is
given by $\f2\pi\sin^2\t d\t$ used in the Sato-Tate conjecture proved in \cite{BLGHT}.

Main results of \cite{K1, K2} give the criterion of the meromorphy of Euler product
$$
L(s,E,H)=\prod_{p:\,\text{prime}} H_{\alpha_\vp^0(p)}(p^{-s})^{-1}
$$
for each polynomial
$$
H(T)\in1+TR(SU(2))[T]
$$
when specialized to the present situation, where $R(SU(2))$ is the ring of virtual characters of $SU(2)$ and
$$
H_{\alpha_\vp^0(p)}(T)\in1+T\C[T]
$$
denotes the polynomial obtained by taking values of coefficients.

We first report Theorems 1 and 2.
The important assumption made in \cite{K1, K2} for
$E=(P(\Z),SU(2),\alpha_\vp^0)$ on the analytic properties of the symmetric power $L$-functions
follow from recently proved automorphy due to Newton-Thorne \cite[Theorem A]{NT1}.
Especially the boundedness in vertical strips of the symmetric power $L$-functions needed in
\cite{K1, K2} follow from Gelbart-Shahidi \cite{GS} (see Shahidi \cite{Sha} and Cogdell-Piatetskii-Shapiro \cite{Co}).

\begin{thm}
Let $H(T)\in1+TR(SU(2))[T]$ be a polynomial of degree $n$.
Then the following properties are equivalent.

\begin{enumerate}[\rm(1)]
\item \vskip -5mm
$H(T)$ is unitary in the sense that there exist functions $\vp_j:\ \R\to\R$ $(j=1,2,...,n)$ such that
$$
H_\t(T)=\prod_{j=1}^n(1-e^{i\vp_j(\t)}T)
$$
for $\t\in[0,\pi]=\mathrm{Conj}(SU(2))$.
\item $L(s,E,H)$ is meromorphic on $\C$.
\end{enumerate}
\end{thm}

\begin{thm}
Let $H(T)\in1+TR(SU(2))[T]$.
Then the following properties are equivalent.
\begin{enumerate}[\rm(1)]
\item \vskip -4mm
$H(T)$ is not unitary.
\item $L(s,E,H)$ is meromorphic in $\Re(s)>0$ with the natural boundary.
\end{enumerate}
\end{thm}

Here we prove a concrete example of degree two for applications in this paper.

\begin{thm}
Let $H(T)=1\pm hT+T^2\in1+TR(SU(2))[T]$.
Then the following properties are equivalent.
\begin{enumerate}[\rm(1)]
\item \vskip -4mm
$H(T)$ is unitary.
\item $|h(\t)|\le2$ for all $\t\in[0,\pi]$.
\item $L(s,E,H)$ is meromorphic on $\C$.
\end{enumerate}
\end{thm}

\noindent{\it Proof.}
Since $(1)\Longleftrightarrow(3)$ is implied in Theorem 1, it suffices to show 
$(1)\Longleftrightarrow(2)$.

\noindent $(1)\Longrightarrow(2)$:\ From
$$
H_\t(T)=(1-e^{i\vp(\t)}T)(1-e^{-i\vp(\t)}T)
$$
for a function $\vp:\ \R\to\R$, we have
$$
|h(\t)|=|2\cos(\vp(\t))|\le2
$$
for $\t\in[0,\pi]$.

\noindent $(2)\Longrightarrow(1)$:\ From $|h(\t)|\le2$ we have
$$
H_\t(T)=\l(1\pm\f{h(\t)+i\sqrt{4-h(\t)^2}}2 T\r)\l(1\pm\f{h(\t)-i\sqrt{4-h(\t)^2}}2 T\r)
$$
and
$$
\l|\f{h(\t)\pm i\sqrt{4-h(\t)^2}}2\r|=1.
$$
\hfill\qed

\section{Proofs of Theorems B and C}
From Theorems 1, 2 and 3 we see that Theorems B and C are proved by the following result.

\begin{thm}
For a monic polynomial $f(x)\in\Z[x]$ of degree $m\ge1$, the following conditions are equivalent.
\begin{enumerate}[\rm(1)]
\item \vskip -4mm
$|f(x)|\le2$ for $-2\le x\le 2$.
\item $f(x)=2T_m(\f x2)$.
\end{enumerate}
\end{thm}

Actually, let
$$
H(T)=1\pm f(2\cos\t)T+T^2\in 1+TR(SU(2))[T],
$$
then we see the following equivalence:
\begin{align*}
\lefteqn{L(s,E,H)\text{ is meromorphic on }\C}\\
&\overset{\text{Th } 1}{\Longleftrightarrow}\quad H(T)\text{ is unitary}\\
&\overset{\text{Th } 3}{\Longleftrightarrow}\quad |f(2\cos\t)|\le2 \text{ for all }\t\in[0,\pi]\\
&\Longleftrightarrow\quad |f(x)|\le2 \text{ for }-1\le x\le 1\\
&\overset{\text{Th } 4}{\Longleftrightarrow}\quad f(x)=2T_m(\f x2)\\
&\Longleftrightarrow\quad L(s,E,H)=\prod_p(1\pm2\cos(m\t(p))p^{-s}+p^{-2s})^{-1}.
\end{align*}
We remark that this also gives a second proof of Theorem A.

\bigskip
\noindent{\it Proof of Theorem 4.}\ 
Let 
$$
g(x)=\f1{2^m}f(2x)\in\R[x].
$$
This is a monic polynomial of degree $m$.
Then, (1) is equivalent to the following condition:
$$
|g(x)|\le\f1{2^{m-1}}\text{ for }-1\le x\le 1.
$$
Such a monic polynomial $g(x)$ is uniquely determined as
$\f1{2^{m-1}}T_m(x)$ by the famous theorem of Chebyshev \cite{C}; see Serre \cite{Se1}.
Hence $f(x)=2T_m\l(\f x2\r).$
\hfill\qed


\begin{bibdiv} \begin{biblist}

\bib{BLGHT}{article}{
   author={T. Barnet-Lamb},
   author={D. Geraghty},
   author={M. Harris},
   author={R. Taylor},
   title={A family of Calabi-Yau varieties and potential automorphy (II)},
   journal={Publ. Res. Inst. Math. Sci},
   volume={47},
   date={2011},
   pages={29-98},
}
\bib{C}{article}{
   author={P.L. Chebyshev},
   title={Th\'eorie des m\'ecanismes connus sous le nom de parall\'elogrammes},
   journal={M\'em. Acad. Sci. P\'etersb.},
   volume={7},
   date={1854},
   pages={539-568},
   note={(= Oe I, 111-143)},
}
\bib{Co}{article}{
   author={Cogdell, J. W.},
   author={Piatetski-Shapiro, I. I.},
   title={Converse theorems, functoriality, and applications to number theory},
   journal={Proceedings of the International Congress of Mathematicians, Volume II},
   date={2002},
   pages={119-128},
}
\bib{D}{article}{
   author={P. Deligne},
   title={La conjecture de Weil (I)},
   journal={Publ. Math. IHES},
   volume={43},
   date={1974},
   pages={273-307},
}
\bib{GS}{article}{
   author={S. Gelbart},
   author={F. Shahidi},
   title={Boundedness of automorphic $L$-functions in vertical strips},
   journal={J. Am. Math. Soc.},
   volume={14},
   date={2001},
   pages={79-107},
}
\bib{K1}{article}{
   author={Kurokawa, N.},
   title={On the meromorphy of Euler products (I)},
   journal={Proc. London Math. Soc. (3)},
   volume={53},
   date={1986},
   pages={1-47},
}
\bib{K2}{article}{
   author={Kurokawa, N.},
   title={On the meromorphy of Euler products (II)},
   journal={Proc. London Math. Soc. (3)},
   volume={53},
   date={1986},
   pages={209-236},
}
\bib{KK}{article}{
   author={Koyama, S.},
   author={Kurokawa, N.},
   title={Rigidity of Euler products},
   journal={arXiv: 2013.06464 v1 [math. NT]  11 Mar},
   date={2021},
}
\bib{L}{article}{
   author={R. P. Langlands},
   title={Problems in the theory of automorphic forms},
   journal={Springer Lecture Notes in Math.},
   volume={170},
   date={1970},
   pages={18-61},
}
\bib{NT1}{article}{
   author={J. Newton},
   author={J. A. Thorne},
   title={Symmetric power functoriality for holomorphic modular forms},
   journal={arXiv: 1912.11261 v2 [math. NT] 17 July},
   date={2020},
}
\bib{NT2}{article}{
   author={J. Newton},
   author={J. A. Thorne},
   title={Symmetric power functoriality for holomorphic modular forms (II)},
   journal={arXiv: 2009.07180 v1 [math. NT] 15 Sep},
   date={2020},
}
\bib{R}{article}{
   author={S. Ramanujan},
   title={On certain arithmetical functions},
   journal={Trans. Cambridge Philosophical Society},
   volume={22},
   date={1916},
   pages={159-184},
}
\bib{Se}{article}{
   author={J. P. Serre},
   title={Une interpr\'etation des congruences relatives \`a la fonction $\tau$ de Ramanujan},
   journal={S\'eminaire Delange-Pisot-Poitou: 1967/68, Th\'eorie des Nombres, Fasc. 1, Exp. 14},
   publisher={Secr\'etariat math\'ematique},
   date={1969},
   pages={17 pp.},
}
\bib{Se1}{article}{
   author={J. P. Serre},
   title={Distribution asymptotique des Valeurs Propres des Endomorphismes de Frobenius [d'apr\`es Abel, Chebyshev, Robinson,…]},
   journal={Ast\'erisque},
   note={S\'eminaire Bourbaki. 2017/2018},
   volume={414},
   date={2019},
   pages={1136-1150},
}
\bib{Sha}{article}{
   author={Shahidi, F.},
   title={Automorphic $L$-functions and functoriality},
   journal={Proceedings of the International Congress of Mathematicians, Volume II},
   date={2002},
   pages={655-666},
}
\bib{S}{article}{
   author={G. Shimura},
   title={On the holomorphy of certain Dirichlet series},
   journal={Proc. London Math. Soc. (3)},
   volume={31},
   date={1975},
   pages={79-98},
}
\bib{T}{article}{
   author={J. Tate},
   title={Algebraic cycles and poles of zeta functions},
   journal={Arithmetical Algebraic Geometry (Proc. Conf. Purdue Univ., 1963)},
   date={1965},
   pages={93-110},
}
\end{biblist} \end{bibdiv}
\end{document}